\documentclass[12pt]{article}
\usepackage{amssymb,amsfonts,amsmath,amsthm, breqn, color, float, mathtools, url, cite, xcolor}
\usepackage{caption}
\newcommand{\witi}{\widetilde}

\newcommand{\cc}{{\mathbb C}}

\newcommand{\beq}{\begin{eqnarray*}}
	\newcommand{\feq}{\end{eqnarray*}}
\newcommand{\beqn}{\begin{eqnarray}}
\newcommand{\feqn}{\end{eqnarray}}

\newtheorem{theorem}{Theorem}
\newtheorem*{conj*}{Conjecture}
\makeatletter \@addtoreset{theorem}{section}\makeatother

\newcommand{\nn}{{\mathbb N}}
\makeatletter \@addtoreset{theorem}{section}\makeatother

\newtheorem{lemma}[theorem]{Lemma}

\newtheorem*{theorema*}{Theorem~A}
\newtheorem*{theoremb*}{Theorem~B}
\newtheorem*{cld*}{Condition $\mbox{LD}_d$}
\newtheorem*{theorem*}{Theorem}

\usepackage{algorithm}
\usepackage{graphicx}
\usepackage[noend]{algpseudocode}
\makeatletter
\def\BState{\State\hskip-\ALG@thistlm}
\makeatother

\DeclareCaptionLabelFormat{noname}{#2}
\newlength\myindent
\def\PP{\mathcal{P}}

\title{Height of records in partitions of a set}

\author{Toufik~Mansour\thanks{ Department of Mathematics, University of Haifa, 199 Abba Khoushy Ave, 3498838 Haifa, Israel;
\newline e-mail: tmansour@univ.haifa.ac.il}
\and
Reza~Rastegar\thanks{Occidental Petroleum Corporation, Houston, TX 77046 and Departments of Mathematics
and Engineering, University of Tulsa, OK 74104, USA - Adjunct Professor; e-mail:  reza\_rastegar2@oxy.com}
\and
Alexander~Roitershtein \thanks{Department of Statistics, Texas A\&M University, College Station, TX 77843, USA;
\newline e-mail: alexander@stat.tamu.edu}
}
\date{August 1, 2019}
\begin{document}
\maketitle

\begin{abstract} 
We study the restricted growth function associated with set partitions, and obtain exact formulas for the number of strong records with height one, the total of record heights over set of partitions, and the number of partitions with a given maximal height of strong records. We also extend some of these results to weak records.
\end{abstract}
	

\section{Introduction and statement of results}
\label{intro} 
We study \textit{records}, also known in the literature as left to right maxima. For a sequence $\pi= \pi_1\cdots \pi_n$ of members of a totally ordered set, $\pi_i$ is a \textit{strong record} if $\pi_i>\max_{1\leq j \leq i-1} \pi_j$. We call it a \textit{weak record} if $\pi_i\geq \max_{1\leq j \leq i-1} \pi_j.$ The interest in records is partially motivated by real-world applications, such as extreme weather studies, tests of randomness, determination of minimal failure, and stresses of electronic components, to name a few. The number of strong records was first studied by R\'enyi \cite{R} for permutations  (see also \cite{Gl,K}) where he proved, among other results, that the number of permutations of $[n]$ with $r$ strong records is equal to the number of such permutations with $r$ $k$-cycles, the latter being given by the unsigned Stirling number of the first kind.
\par 
More recently, Myers and Wilf \cite{MW} extended the study of records to multiset permutations and words. They derived the generating function for the number of permutations of a fixed multiset which contains exactly $r$ (strong or weak) records. They also gave the average number of (strong and weak) records among all permutations of the multiset. Their work has been extended to different totally ordered discrete sets and a variety of statistics of records. For example, Kortchemski in \cite{Ko} studied the asymptotic behavior of the number of permutations having $r$ records and the permutations for which the sum of the positions of their records is $r$. Knopfmacher, Mansour and Wagner \cite{KMW} studied the sum of the positions of records in set partitions. Caki\'c, Mansour and Smith \cite{CMS} investigated further another statistic of records in set partitions. Asakly \cite{Asa} obtained an explicit formula and asymptotic estimate for the total number of sum of weighted records over set partitions of $[n]$ in terms of Bell numbers.
\par
In this paper we study the height of strong and weak records in set partitions. Most of our results are exact and obtained by analyzing underlying generating functions. The asymptotic behavior of the quantities of interest can be subsequently deduced from the exact formulas.
\par
We begin with notations. A \textit{partition} of a set $A$ is a collection of non-empty, mutually disjoint subsets, called \textit{blocks},
whose union is the set $A$. A partition $\Pi$ with $k$ blocks is called a \textit{$k$-partition} and denoted by $\Pi = A_1\mid A_2\mid\ldots\mid A_k$. 
Let $\PP_{n,k}$ be the set of $k$-partitions of $[n]=\{1,2,\ldots,n\}$
and define $\PP_n =\cup_{k=0}^n\PP_{n,k}$ to be the collection of all partitions of $[n]$. It is well-known that 
$\mbox{Card}(\PP_{n,k})=S_{n,k}$ and $\mbox{Card}(\PP_{n})=B_n,$ where $S_{n,k}$ 
is a Stirling number of the second kind  and $B_n$ is the $n$-th Bell number. The Stirling numbers 
of the first kind $s_{n,k}$ and the second kind $S_{n,k}$ can be introduced algebraically in several different ways, including 
formulas \eqref{snkIdent} and \eqref{SnkIdent} given below. Alternatively, one can define the sequence of Stirling numbers of the first kind 
as the solution to the recursion 
\beq
s_{n,k}=-(n-1)s_{n-1,k}+s_{n-1,k-1},\qquad n,k\in\nn,\,k\leq n,
\feq
with the initial conditions $s_{0,0}=1$ and $s_{0,n}=0,$ and the sequence of Stirling numbers of the second kind as the solution to the recursion 
\beq
S_{n,k}=kS_{n-1,k}+S_{n-1,k-1},\qquad n,k\in\nn,\,k\leq n,
\feq
with the initial conditions $S_{0,0}=1$ and $S_{0,n}=0.$ The sequence of Bell numbers $(B_n)_{n\geq 0}$ can be then defined through the formula 
$B_n=\sum_{k=0}^n S_{n,k},$ or, for instance, recursively via the formula $B_{n+1}=\sum_{k=0}^n \binom{n}{k}B_k$ 
with the initial condition $B_0=1.$ For more on the Stirling and Bell numbers see, for instance, \cite{comtet, M, S} and references therein.     
\par 
A $k$-partition $\Pi$ is said to be in the \textit{standard form} if the blocks $A_i$ are labeled in such a way that
\beqn
\label{min}
\min A_1 < \min A_2<\cdots< \min A_k. 
\feqn 
Note that if a $k$-partition of $[n],$ call it $\Pi,$ is in 
the standard form, then it can be represented equivalently by the 
canonical sequential form $\Pi= \pi_1\pi_2\ldots\pi_n$ where 
$\pi_i\in [n]$ and $i\in A_{\pi_i}$ for all $i$ (see, for instance, \cite{M,S}). In words, $\pi_i$ is the label of the 
partition block that contains $i.$ We remark that the canonical form
is sometimes referred as a \textit{restricted growth function} (see, for instance, \cite{M}).
\par Notices that \eqref{min} guarantees that $1\in A_1.$  
For example, the canonical representation of the partition $\{1,4,7\}\mid \{2,3,6,9\}\mid \{5,8\}\in\PP_{9,3}$ is $122132132.$ 
Throughout  this paper we consider the canonical form as a representation of partitions.
Note that a word $\pi\in [k]^n$ is a canonical representation of a $k$-partition of $[n]$ if and only if the following holds:
\begin{enumerate}
\item Each element of $[k]$ appears at least once in $\pi.$
\item For all $1\le i<j\le k,$ the first occurrence of $i$ precedes the first occurrence of $j$.
\end{enumerate}	 
To define \emph{records} for a $k$-partition $\Pi$ of $[n]$, we represent $\Pi$ by its corresponding restricted growth function $\pi=\pi_1\cdots \pi_n$. Then, we define the records of $\Pi$ as records of their restricted growth function $\pi$. For a record, $\pi_i,$ we call $i$ the position of the record $\pi_i.$ Similarly, we define the height of this record as $\pi_i-\pi_{i-1}$ whenever $i>1$ and (within the framework of strong records) $0$ if $i=1$.
\par	 	
Our first result is
\begin{theorem}
\label{main_thm1}
$\mbox{}$
\item [(i)] Let $n\geq k\geq 1$ and $r\geq0$. The number of $k$-partitions of $[n]$ with exactly $r$ records of height one is
\begin{align*}
\sum_{j=0}^{k-1}(-1)^{r-j}\binom{k-1-j}{r-j}|s_{k-1,j}|S_{n-k+1+j,k}.
\end{align*}
\item [(ii)] 
The total number of strong records of height one in all partitions in $\PP_{n,k}$ is given by
\beq 
\frac{1}{2}S_{n+1,k}+\frac{1}{2}S_{n,k}-S_{n-1,k}-S_{n-1,k-1}-\frac{1}{2}S_{n-1,k-2}.
\feq
\item [(iii)] 
The total number of strong records of height one in all partitions in $\PP_{n}$ is given by
\beqn
\label{as}
\frac{1}{2}B_{n+1}+\frac{1}{2}B_n-\frac{5}{2}B_{n-1} = \left(\frac{n+1}{2\xi_n}+\frac{1}{2}+ O\Bigl(\frac{\log n}{n}\Bigr)\right)B_n,
\feqn
where $\xi_n$ is the unique positive root of the equation $\xi e^{\xi}=n+1.$
\end{theorem}
We remark that (see, for instance, Example~VIII.6 in \cite{flaj})
\beq
\xi_n = \log n - \log \log n + \frac{\log \log n}{\log n} + O\bigg(\frac{\log^2 \log n}{\log^2 n}\bigg).
\feq
Similarly to Theorem~\ref{main_thm1}, next theorem provides some insight on strong records with the exception that we drop the assumption of height one. 	
\begin{theorem}
\label{main_thm2}
$\mbox{}$
\item [(i)] The sum of the heights of the strong records over all the partitions in $\PP_{n,k}$ is given by
\beq
(k-1)S_{n,k}+\binom{k}{3}S_{n-1,k}.
\feq
\item [(ii)] The sum of the heights of the strong records over all the partitions in  $\PP_n$
is given by
\beqn
\label{as1}
\frac{1}{6}B_{n+2}-\frac{2}{3}B_{n}-\frac{1}{6}B_{n-1}
=
\bigg(\frac{(n+2)(n+1)}{6\xi^2}-\frac{2}{3}+ O\Bigl(\frac{\log n}{n}\Bigr) \bigg)B_n.
\feqn 
\item [(iii)] Let $\PP^h_{n,k}$ be the set of partitions in $\PP_{n,k}$ whose maximum height is $h.$ Then	
\beq
\mbox{Card}(\PP^h_{n,k}) = \sum_{j=1}^{k-h}s_{k-h,j}S_{n-(k-h)+j,k}.
\feq
\end{theorem}
Next theorem extends our previous results to weak records.	
\begin{theorem}
\label{main_thm3}
$\mbox{}$
\item [(i)]  Let $n\geq k\geq 1$ and $r\geq0$. The number of $k$-partitions of 
$[n]$ with exactly $r$ weak records of height one is given by
\beq
&&\sum_{i=0}^{k-1}\sum_{a=0}^{(n+i-k)/2}(-1)^{m+a-i}\binom{n+i-a-k}{a}\binom{k+a-1-i}{m-i}S_{n+1+i-2a-k,k}|s_{k-1,i}|\\
&&\qquad +\sum_{i=0}^{k-1}\sum_{a=0}^{(n+i-1-k)/2}\sum_{j=0}^{n+i-1-2a-k}(-1)^{m+a-i}\binom{j-1+a}{a}\binom{k+a-i}{m-i}S_{j,k}|s_{k-1,i}|.
\feq
\item [(ii)] The total number of weak records of height one in all the partitions in $\PP_{n,k}$ is given by 
\beq 
\frac{1}{2}S_{n+2,k}+\frac{1}{2}S_{n+1,k}+nS_{n-1,k}-S_{n,k-1}
-\frac{1}{2}S_{n,k-2}-S_{n-1,k}-\sum_{j=0}^{n-1}\binom{n}{j}S_{n-1-j,k-1}.
\feq
\item [(iii)] 
The total number of weak records of height one in all 
the partitions in $\PP_n$ is given by 
\beq &&\frac{1}{2}B_{n+1}+\frac{1}{2}B_{n}+(n-1)B_{n-2}-\frac{3}{2}B_{n-1}-B_{n-2}
-\sum_{j=0}^{n-1}\binom{n}{j+1}B_j 
\\ 
&& \qquad \qquad = \bigg(\frac{n+1}{2\xi_n}+\frac{1}{2}+
O\Bigl(\frac{\log n}{n}\Bigr)\bigg)B_{n}.
\feq
\item [(iv)] The sum of the heights of the weak records over all the partitions in $\PP_{n,k}$ is given by 
\beq
(k-1)S_{n,k}+\binom{k}{3}S_{n-1,k}+\binom{k+1}{3}S_{n-2,k}.
\feq
\item [(v)] The sum of the heights of the weak records over all the partitions in $\PP_n$ is given by
\beq &&\frac{1}{6}B_{n+2}+\frac{1}{6}B_{n+1}-\frac{7}{6}B_{n}-\frac{1}{3}B_{n-1}+\frac{1}{3}B_{n-2}\\
&& \qquad = \bigg(\frac{(n+2)(n+1)}{6\xi_n^2}+\frac{(n+1)}{6\xi_n}-\frac{7}{6} + O\Big(\frac{\log n}{n}\Big)\bigg)B_{n}.
\feq
\end{theorem} 
The rest of the paper is organized as follows. Theorems~\ref{main_thm1} and~\ref{main_thm2} 
are proved in Section~\ref{strong_rec_sect}, and the proof of Theorem \ref{main_thm3} is given in 
Section~\ref{weak_rec_sect}. The proofs are analytical and rely on the analysis of generating functions.  
	
\section{Strong records. Proof of Theorems~\ref{main_thm1} and~\ref{main_thm2}}
\label{strong_rec_sect}
In this section we are concerned with the height of strong records. The section is divided into four subsection. 
Section~\ref{prelims} is of a preliminary function, its goal is to introduce certain generating function associated with strong records. 
The proof of Theorem~\ref{main_thm1}  is given in Section~\ref{sr1}. 
The proof of claims (i) and (ii) of Theorem~\ref{main_thm2} is contained in Section~\ref{totals}, and 
the proof of part of the theorem is included in Section~\ref{maxs}. 
\subsection{Preliminaries} 
\label{prelims}
We first introduce a few notations.
Throughout this paper, for any given ordinary generating function
\beq
A(x)=\sum_{n\geq0} a_n x^n,\qquad x\in \cc,
\feq
we use $\witi A$ to denote the corresponding exponential generating function. That is,
\beqn
\label{Ae}
\witi A(x):=\sum_{n\geq0}a_n\frac{x^n}{n!}.
\feqn
We define $P_{n,k}=P_{n,k}(q_1,q_2,\ldots,q_{k-1})$ to be the generating function of the number of
partitions in $\PP_{n,k}$ with respect to the height of the strong records, where $q_i$ counts the number of records of the height $i$.
For $k\geq 0$ define
\begin{align*}
P_k(x)=P_k(x,q_1,\ldots,q_{k-1})=\sum_{n\geq k}P_{n,k}x^n.
\end{align*}
We remark that $P_0(x)=1$ and $P_1(x)=\frac{x}{1-x}$.
Thus, unless explicitly said otherwise, for the rest of this paper we assume that $k\geq2$.
\par
In order to study the generating function $P_k(x)$ we define
\begin{align*}
P_k(x|a_sa_{s-1}\cdots a_1)=\sum_{n\geq k}P_{n,k}(q_1,\ldots,q_{k-1}|a_sa_{s-1}\cdots a_1)x^n,
\end{align*}
where
\begin{align*}
P_{n,k}(q_1,\ldots,q_{k-1}|a_sa_{s-1}\cdots a_1)
\end{align*}
is the generating function for the number of $k$-partitions of the form $\pi=\pi'a_sa_{s-1}\cdots a_1$ in $\PP_{n,k}$.
Observe that the last symbol in the canonical representation of a $k$-partition  cannot be a strong record unless it equals
to $k.$ Hence, for $a\in [k-1]$ we have
\begin{align}\label{eqA1}
P_k(x|a)=xP_k(x).
\end{align}
Using the fact that $P_k(x)=\sum_{a=1}^kP_k(x|a),$ we obtain the identity
\begin{align}
\label{eqA2}
P_k(x)=(k-1)xP_k(x)+P_k(x|k).
\end{align}
Furthermore,
\beq
P_k(x|k)&=&xP_k(x)+x\sum_{i=1}^{k-1} q_iP_{k-1}(x|k-i) \\
&=& xP_k(x)+q_1P_{k-1}(x|k-1) + x^2P_{k-1}(x) \sum_{i=2}^{k-1} q_i,
\feq
where the last equality follows from \eqref{eqA1}, and the empty sum $\sum_{i=2}^{k-1} q_i$ is considered to be zero when $k=2.$ 
It follows then from \eqref{eqA2} that
\beq
(1-(k-1)x)P_k(x)=xP_k(x)+xq_1(1-(k-2)x)P_{k-1}(x) && \\ \qquad + x^2(q_2+\cdots+q_{k-1})P_{k-1}(x).
\feq
Equivalently, 
\beqn
\nonumber  
P_k(x)&=&x\frac{q_1+x\bigl(q_2+\cdots+q_{k-1}-(k-2)q_1\bigr)}{1-kx}P_{k-1}(x)
\\
\label{Plx}
&=&
\ldots=	
x^k\frac{\prod_{j=2}^k(q_1+x(q_2+\cdots+q_{j-1}-(j-2)q_1))}{\prod_{j=1}^k(1-jx)}.
\feqn
We remark that if $q_1= \cdots = q_{k-1}=q$, then \eqref{Plx} becomes
\beq
P_k(x;q,q,\ldots,q)=\frac{x^kq^{k-1}}{\prod_{j=1}^k(1-jx)}=
q^{k-1}\sum_{n\geq0}S_{n,k}x^n,
\feq
where the second equality is due to \eqref{SnkIdent}. 
This is expected because $\mbox{Card}(\PP_{n,k})=S_{n,k}$ and any set $k$-partition of $[n]$ has exactly $k-1$ strong records.	
\subsection{Strong records of height one. Proof of Theorem~\ref{main_thm1}}
\label{sr1}
The goal os this section is to prove Theorem~\ref{main_thm1}. 
In order to analyze the number of strong records of height one, we let $q_1=q$ and $q_i=1$ for all $i\geq2$ in \eqref{Plx}.
Denote
\beq
T_k(x;q):=P_k(x;q,1,\ldots,1).
\feq
Then \eqref{Plx} yields the recursion
\beqn
\label{recTxq}
(1-kx)T_k(x;q)=x\big(q+(k-2)x(1-q)\big)T_{k-1}(x;q)
\feqn
with $T_0(x;q)=1$ and $T_1(x;q)=\frac{x}{1-x}$. Thus,
\beqn
T_k(x;q)&=&\frac{x^k\prod_{j=0}^{k-2}(q+jx(1-q))}{\prod_{j=0}^k(1-jx)}
\nonumber
\\
&=&
\frac{x^{2k-1}(1-q)^{k-1}\prod_{j=0}^{k-2}\left(\frac{q}{x(1-q)}+j\right)}{\prod_{j=0}^k(1-jx)}.
\label{Tkxq}
\feqn
Recall (see, for instance, \cite{comtet, M}) that Stirling numbers satisfy the identities
\beqn
\label{snkIdent}
\prod_{j=0}^{n-1}(x+j)=\sum_{j=0}^n|s_{n,j}|x^j\qquad \mbox{\rm and}\qquad \prod_{j=1}^k (1-jy)=\sum_{j=0}^ks_{k+1,j+1}y^{k-j}.
\feqn
(note that the latter identity is equivalent to the former one with $x=-1/y$), and
\beqn
\label{SnkIdent}
\frac{x^k}{\prod_{j=1}^k(1-jx)}=\sum_{n\geq0}S_{n,k}x^n.
\feqn
Therefore, \eqref{Tkxq} can be expressed as
\begin{align*}
T_k(x;q)&=\sum_{j=0}^{k-1}|s_{k-1,j}|q^jx^{k-1-j}(1-q)^{k-1-j}
\sum_{n\geq k}S_{n,k}x^n
\\ 
&=\sum_{n\geq k}\sum_{j=0}^{k-1}\sum_{i=0}^{k-1-j}(-1)^i
\binom{k-1-j}{i}|s_{k-1,j}|S_{n,k}q^{i+j}x^{n+k-1-j}.
\end{align*}
Changing the order of summation from $\sum_{n\geq k}\sum_{j=0}^{k-1}\sum_{i=0}^{k-1-j}...$ to
\beq
\sum_{r=0}^{k-1}\bigg\{\sum_{N=k-1-r}^{k-1}\sum_{j=0}^{k-1-N}...+\sum_{N=k}^\infty\sum_{j=0}^r...\bigg\}
\feq
and inspecting the coefficient of $q^rx^n$ in the resulting formula for $T_k(x;q)$ finishes the proof of Theorem \ref{main_thm1}-(i).	
\par 	
Next, we prove Theorem \ref{main_thm1}-(ii) and ~(iii). Let $\witi T$ be the exponential generating function associated with $T,$ 
as defined in \eqref{Ae}. Then, by \eqref{recTxq}, we get
\beqn
\frac{\partial^2}{\partial x^2}\witi T_k(x;q) -k\frac{\partial}{\partial x}\witi T_k(x;q)
=q\frac{\partial}{\partial x}\witi T_{k-1}(x;q)+(k-2)(1-q)\witi T_{k-1}(x;q)
\feqn
with $\witi T_0(x;q)=1$ and $\witi T_1(x;q)=e^x-1$. 
Define 
\beq
\witi T(x,y;q)=\sum_{k\geq0}\witi T_k(x;q)y^k.
\feq 
Multiplying by $y^k$ and summing up over $k\geq2$ gives
\beqn
&& \frac{\partial^2}{\partial x^2}\witi T(x,y;q) -y\frac{\partial^2}{\partial x\partial y}\witi T(x,y;q) 
\nonumber \\
&& \quad =qy\frac{\partial}{\partial x}\witi T(x,y;q)-(1-q)y(\witi T(x,y;q)-1)+(1-q)y^2\frac{\partial}{\partial y}\witi T(x,y;q). \label{eqh1a}
\feqn 
Note that $\witi T(x,y;1)=e^{y(e^x-1)}.$ By taking the derivative in the both sides of \eqref{eqh1a} 
and evaluating the result at $q=1,$ we obtain the equation
\beq
\Bigl\{\frac{\partial^3}{\partial x^2\partial q}\witi T(x,y;q) -y
\frac{\partial^3}{\partial x\partial y\partial q}\witi T(x,y;q)-y\frac{\partial^2}{\partial x\partial q}\witi T(x,y;q)\Bigr\}\,\Big|_{q=1}=y(y+1)e^{y(e^x-1)}-y,
\feq
with the initial condition $$\frac{\partial^2}{\partial x\partial q}\witi T(x,y;q)\,\Bigm|_{q=1}=0.$$ Therefore,
\beq 
\frac{\partial^2}{\partial x\partial q}\witi T(x,y;q)\,\Bigm|_{q=1} 
=\frac{1}{2}\frac{\partial^2}{\partial x^2}e^{y(e^x-1)}+
\frac{1}{2}\frac{\partial}{\partial x}e^{y(e^x-1)}-(1+y+y^2/2)e^{y(e^x-1)}+1.
\feq 
Inspecting the coefficient of $x^ny^k/n!$ in $\frac{\partial}{\partial q}\witi T(x,y;q)\mid_{q=1}$ 
and the coefficient of $x^n/n!$ in $\frac{\partial}{\partial q}\witi T(x,1;q)\mid_{q=1}$ 
finishes the proof of Theorems~\ref{main_thm1}-(ii) and~(iii).
\par
We remark that in order to derive the asymptotic estimate in Theorem~\ref{main_thm1}-(iii), we use the fact that (see, for instance, \cite{Can})
the identity
\begin{align}
\label{asym}
B_{n+h}=B_n\frac{(n+h)!}{n!\xi_n^h}\Bigl(1+O\Bigl(\frac{\log n}{n}\Bigr)\Bigr)
\end{align}
holds uniformly for $h=O(\log n)$, where $\xi_n$ is the unique positive root of $\xi e^\xi=n+1$. \hfill\hfill\qed
	
\subsection{Total height of strong records. Proof of Theorem~\ref{main_thm2}-(i),(ii)}
\label{totals}
The goal of this section is to prove part (i) and part (ii) of Theorem~\ref{main_thm2}. 
In order to analyze the sum of heights of all strong records over $P_{k,n}$, we let $q_i=q^i$ in \eqref{Plx}.
Define $T_k(x;q):=P_k(x;q,q^2,\ldots,q^{k-1})$ to be the generating function for the number of partitions of 
$\PP_{n,k}$ accounting for the total heights of all strong records. In this notation, \eqref{Plx} becomes 
\beqn 
\label{Tkxq2}
T_k(x;q)=x^kq^{k-1}\frac{\prod_{j=2}^k(1+x(q+\cdots+q^{j-2}-(j-2)))}{\prod_{j=1}^k(1-jx)}.
\feqn
Take the derivative with respect to $q$ on the both sides of \eqref{Tkxq2} and evaluate the result at $q=1$. Then
\beq 
\frac{\partial}{\partial q}T_k(x;q)\Big|_{q=1} 
=\frac{(k-1)x^k}{\prod_{j=1}^k(1-jx)}+\frac{x^{k+1}\binom{k}{3}}{\prod_{j=1}^k(1-jx)}.
\feq
Recall \eqref{SnkIdent}. An argument similar to the one given in the previous section finishes the proof of Theorem \ref{main_thm2}-(i).	
\par 
Before we proceed with the proof of part (ii) of Theorem \ref{main_thm2}, 
we recall another identity for $S_{n,k}$. We refer the reader to \cite{M} for a review of this type of identities 
for Stirling's numbers. Recall $\frac{(e^x-1)^k}{k!}=\sum_{n\geq k} S_{n,k}\frac{x^n}{n!}.$
Then
\begin{align*}
&\sum_{k\geq2}\sum_{n\geq k}\Big((k-1)S_{n,k}+\binom{k}{3}S_{n-1,k}\Big)\frac{x^ny^k}{n!}
\\
&
\qquad 
=\sum_{k\geq2}\left(\frac{k-1}{k!}(e^x-1)^k+\binom{k}{3}\int_0^x\frac{(e^t-1)^k}{k!}dt\right)y^k
\\
&
\qquad 
=1+(y(e^x-1)-1)e^{y(e^x-1)}+\int_0^x\frac{y^3}{6}(e^t-1)^3e^{y(e^t-1)}dt.
\end{align*}
Define $\witi T_k(x):=\frac{\partial}{\partial q}T_k(x;q)\big|_{q=1}$. Then, \eqref{Tkxq2} implies that
\begin{align*}
\frac{\partial}{\partial x}\sum_{k\geq0}\witi T_k(x)y^k
&=y^2(e^x-1)e^{y(e^x+x-1)}+\frac{y^3}{6}(e^x-1)^3e^{y(e^x-1)}.
\end{align*}
In particular, 
\begin{align*} 
\frac{\partial}{\partial x}\sum_{k\geq0}\witi T_k(x)
&=(e^x-1)e^{e^x+x-1}+\frac{1}{6}(e^x-1)^3e^{e^x-1}
\\
&=\frac{1}{6}\frac{d^3}{dx^3}e^{e^x-1}-\frac{2}{3}
\frac{d}{dx}e^{e^x-1}-\frac{1}{6}e^{e^x-1}.
\end{align*} 
The inspection of the power series expansion of the above expression for $\frac{\partial}{\partial x}\sum_{k\geq0}\witi T_k(x)$ reveals 
that the coefficient of $\frac{x^n}{n!}$ is given by
\beq 
\frac{1}{6}B_{n+3}-\frac{2}{3}B_{n+1}-\frac{1}{6}B_n.
\feq
This completes the proof of Theorem~\ref{main_thm2}-(ii). 
Note that the asymptotic estimate in \eqref{as1} is a straightforward consequence of the exact formula and \eqref{asym}.\hfill\hfill\qed
\subsection{Maximum height of strong records. Proof of Theorem~\ref{main_thm2}-(iii)}
\label{maxs}
We now turn to the proof of Theorem~\ref{main_thm2}-(iii). Let $\PP^h_{n,k}$ 
be the set of partitions in $\PP_{n,k}$ whose maximum height is $h$. 
To obtain the generating function of $\mbox{Card}(\PP^h_{n,k}),$ 
we let $q_i=1$ for $i\in[h]$ and $q_i=0$ for $i\geq h+1$ in \eqref{Plx},
thus setting 
\beq 
M_k(x):=P_k(x;q_1,q_2,\ldots,q_{k-1})=
x^k\frac{\prod_{j=2}^k(q_1+x(q_2+\cdots+q_{j-1}-(j-2)q_1))}{\prod_{j=1}^k(1-jx)}.
\feq
Hence,
\beq
M_k(x)=
\left\{
\begin{array}{ll}
\frac{x^k}{\prod_{j=1}^k(1-jx)} &\quad\mbox{\rm for}\quad k\leq h+1,
\\
[10pt]		
\frac{x^k(1-x)\cdots(1-(k-h-1)x)}{\prod_{j=1}^k(1-jx)} &\quad\mbox{\rm otherwise.}
\end{array}
\right.
\feq
It follows then from \eqref{snkIdent} and \eqref{SnkIdent} that the number of partitions in $\PP^h_{n,k}$ is given by
\beq 
\sum_{j=1}^{k-h}s_{k-h,j}S_{n-(k-h)+j,k}.
\feq 
This completes the proof of Theorem~\ref{main_thm2}-(iii). \hfill\hfill\qed 
\section{Weak records. Proof of Theorem~\ref{main_thm3}}
\label{weak_rec_sect}
This section is devoted to to prove Theorem \ref{main_thm3}.  The proof of claims (i)--(iii) of the theorem 
is given in Section~\ref{wrh1}. The proof of claims (iv) and (v) is included in Section~\ref{wr}. Section~\ref{prelimg} 
is of a preliminary function, its goal is to compute a certain generating function associated with weak records.    
\subsection{Generating functions}  
\label{prelimg}  
The goal of this section is to obtain Lemma~3.1, a closed form expression for certain generating functions. 
Let $Q_{n,k}=Q_{n,k}(q_1,q_2,\ldots,q_{k-1})$ be the multivariate generating function for the number of partitions in
$\PP_{n,k}$ accounting for the heights of weak records. Here $q_i$ counts the number of the weak records of height $i$.
Define
\beq
Q_k(x)=Q_k(x,q_1,\ldots,q_k)=\sum_{n\geq k}Q_{n,k}x^n.
\feq
It is easy to see that $$Q_0(x)=1\quad \mbox{and} \quad Q_1(x)=\frac{x}{1-x}.$$
From now we will assume that $k\geq2$.
In order to study $Q_k(x)$, we introduce auxiliary generating functions
\beq 
Q_k(x|a_sa_{s-1}\cdots a_1)=\sum_{n\geq k} 
Q_{n,k}(q_1,\ldots,q_{k-1}|a_sa_{s-1}\cdots a_1)x^n,
\feq 
where
$Q_{n,k}(q_1,\ldots,q_{k-1}|a_sa_{s-1}\cdots a_1)x^n$ is the 
generating function for the number of partitions 
$\pi=\pi'a_sa_{s-1}\cdots a_1$ in $\PP_{n,k}$ accounting for 
the heights of the weak records.	 
\par
Let $\pi=\pi'a$ be a partition with $k$ blocks and $1\leq a\leq k-1$. 
Since the partition is assumed to be in the standard form and $a<k,$ it cannot be a record. Thus
\begin{align}
\label{eqB1}
Q_k(x|a)=xQ_k(x).
\end{align}
Therefore,
\begin{align}
\label{eqB2}
Q_k(x)=\sum_{a=1}^kQ_k(x|a)=(k-1)xQ_k(x)+Q_k(x|k).
\end{align}
In order to compute $Q_k(x|k)$, we consider $k$-partitions in the form $\pi=\pi'k$ 
and two separate scenarios: when $\pi'$ is a $k$-partition and when it is a $(k-1)$-partition. We have:
\beq
Q_k(x|k)&=&xq_{k-1}Q_k(x|1)+\cdots+xq_1Q_k(x|k-1)+xQ_k(x|k)
\\
&& \quad +xq_{k-1}Q_{k-1}(x|1)+xq_{k-2}Q_{k-1}(x|2)+\cdots+xq_1Q_{k-1}(x|k-1).
\feq
With the help of \eqref{eqB1}, this formula simplifies to
\beq
Q_k(x|k)&=&xQ_k(x|k)+x^2(q_1+\cdots+q_{k-1})Q_k(x)
\\
&& \quad +xq_1Q_{k-1}(x|k-1)+x^2(q_2+\cdots+q_{k-1})Q_{k-1}(x).
\feq
It then follows from \eqref{eqB2} that
\begin{align*}
Q_k(x)=\frac{xq_1(1-(k-2)x)+x^2(q_2+\cdots+q_{k-1})}{(1-x)(1-(k-1)x)-x^2(q_1+\cdots+q_{k-1})}Q_{k-1}(x).
\end{align*}
Finally, by induction on $k$ we get
\begin{lemma}
\label{thB1}
For $k\geq1$,
\beqn\label{Qkqs} Q_k(x;q_1,q_2,\ldots,q_{k-1})=
x^k\frac{\prod_{j=2}^k\Big(q_1\big(1-(j-2)x\big)+x\big(q_2+\cdots+q_{j-1}\big)\Big)}{\prod_{j=1}^k
\Big(\big(1-x)(1-(j-1)x\big)-x^2\big(q_1+\cdots+q_{j-1}\big)\Big)}.
\feqn
\end{lemma}	
We remark that the particular case of \eqref{Qkqs} with $q_1=\cdots=q_{k-1}=q,$ namely
\beq 
Q_k(x;q,q,\ldots,q)=\frac{q^{k-2}x^k}{\prod_{j=1}^k\big(1-jx+(j-1)x^2(1-q)\big)},
\feq
has been calculated in \cite{KMW}.
\subsection{Weak records of height one. Proof of Theorem~\ref{main_thm3}-(i)--(iii)}
\label{wrh1}		
In order to study the weak records of height one, we will utilize Lemma~\ref{thB1} with $q_1=q$ 
and $q_i=1$ for all $i\geq2$. It follows from \eqref{Qkqs} that
\beq 
Q_k(x;q,1,\ldots,1)&=&\frac{x^k}{1-x}\frac{\prod_{j=2}^k(q+x(j-2)(1-q))}{\prod_{j=2}^k(1-jx+(1-q)x^2)} 
\\ 
&=&
\frac{1-x+(1-q)x^2}{(1-x)}\,
\frac{x^k/\bigl(1+(1-q)x^2\bigr)^k}{\prod_{j=1}^k\big(1-\frac{jx}{1+(1-q)x^2}\big)}
\, 
\prod_{j=0}^{k-2}\big(q+jx(1-q)\big). 
\feq
Next, using \eqref{snkIdent} and \eqref{SnkIdent}, we obtain: 
\begin{align*}
&Q_k(x;q,1,\ldots,1)
\\
&=\frac{1-x+(1-q)x^2}{1-x}\sum_{j\geq k}S_{j,k}\frac{x^j}{(1+(1-q)x^2)^j}
\sum_{i=0}^{k-1}|s_{k-1,i}|q^ix^{k-1-i}(1-q)^{k-1-i}
\\
&=\left(1-\frac{(q-1)x^2}{1-x}\right)\sum_{j\geq k} \sum_{i=0}^{k-1}\sum_{a\geq0}(-1)^{k-1-i}\binom{j-1+a}{a}S_{j,k}|s_{k-1,i}|q^ix^{k+2a+j-1-i}(q-1)^{k+a-1-i}. 
\end{align*}
Thus, the coefficient of $q^r$ in $Q_k(x;q,1,\ldots,1)$ is
\begin{align*}
&\sum_{j\geq k}\sum_{i=0}^{k-1}\sum_{a\geq0}(-1)^{r+a-i}\binom{j-1+a}{a}\binom{k+a-1-i}{r-i}S_{j,k}|s_{k-1,i}|x^{k+2a+j-1-i}\\
&\qquad +\frac{1}{1-x}\sum_{j\geq k}\sum_{i=0}^{k-1}\sum_{a\geq0}(-1)^{r+a-i}\binom{j-1+a}{a}\binom{k+a-i}{r-i}S_{j,k}|s_{k-1,i}|x^{k+2a+j+1-i}.
\end{align*}
Finally, inspecting the coefficient of $x^n$ we obtain Theorem~\ref{main_thm3}-(i).
\par
Next, we prove Theorem~\ref{main_thm3}-(ii) and~(iii). To this end, let $q_1=q$ and $q_i=1$ for all $i\geq2$ in \eqref{thB1}. 
Then
\beq 
Q_k(x;q,1,\ldots,1)=\frac{x^k}{1-x}
\frac{\prod_{j=2}^k(q+x(j-2)(1-q))}{\prod_{j=2}^k(1-jx+(1-q)x^2)}.
\feq
This can be written recursively as
\beqn 
\label{recQk}
(1-kx+(1-q)x^2)Q_k(x;q,1,\ldots,1)=x(q+x(k-2)(1-q))Q_{k-1}(x;q,1,\ldots,1).
\feqn
Let $\witi Q$ be the exponential generating function associated with $Q,$ as defined by \eqref{Ae}. 
Then \eqref{recQk} gives the following integral equation:
\begin{align*}
\witi Q_k(x;q)&-k\int_0^x\witi Q_k(t;q)dt+
(1-q)\int_0^x\int_0^t\witi Q_k(r;q)dr\,dt
\\
&=q\int_0^x\witi Q_{k-1}(t;q)dt+(k-2)(1-q)\int_0^x\int_0^t\witi Q_{k-1}(r;q)dr \,dt.
\end{align*}
Differentiating twice with respect to $x,$ we arrive to the following second order ODE:
\beqn
&& \frac{\partial^2}{\partial x^2}\witi Q_k(x;q)-k
\frac{\partial}{\partial x}\witi Q_k(x;q)+(1-q)\witi Q_k(x;q) 
\nonumber 
\\
&&\quad = q\frac{\partial}{\partial x}\witi Q_{k-1}(x;q)+(k-2)(1-q)\witi Q_{k-1}(x;q) 
\label{QtildeODE}
\feqn
with initial conditions 
\beq
\witi Q_0(x;q)=1 \quad \text{and} \quad \witi Q_1(x;q)=e^x-1.
\feq
Define $\witi Q(x,y;q)=\sum_{k\geq0}\witi Q_k(x;q)y^k$. 
Multiplying both sides of \eqref{QtildeODE} by $y^k$ and summing over $k\geq2$, we get
\begin{align}
\frac{\partial^2}{\partial x^2}\witi Q(x,y;q)-y\frac{\partial^2}{\partial x\partial y}\witi Q(x,y;q)&+
(1-q)((1+y)\witi Q(x,y;q)-1-e^xy)
\nonumber
\\
&=qy\frac{\partial}{\partial x}\witi Q(x,y;q)+(1-q)y^2\frac{\partial}{\partial y}\witi Q(x,y;q).\label{eqBh1}
\end{align}
Let $\witi T(x,y)=\frac{\partial}{\partial q}\witi Q(x,y;q)\mid_{q=1}$. 
Since $\witi Q(x,y;1)=e^{y(e^x-1)}$, \eqref{eqBh1} gives
\begin{align*}
\frac{\partial^2}{\partial x^2}\witi T(x,y)-y\frac{\partial^2}
{\partial x\partial y}\witi T(x,y)-y\frac{\partial}{\partial x}\witi T(x,y)=(1+y+y^2)e^{y(e^x-1)}-1-e^xy.
\end{align*}
Finally, using the initial condition $\frac{\partial}{\partial x}\witi T(0,y)=0$ we arrive to
\beq
\frac{\partial}{\partial x}\witi T(x,y)&=&\frac{1}{2}\frac{\partial^2}{\partial x^2}e^{y(e^x-1)}
+\frac{1}{2}\frac{\partial}{\partial x}e^{y(e^x-1)}
\\
&& \quad +(x-y-y^2/2)e^{y(e^x-1)}-(1+ye^x)\int_0^xe^{y(e^t-1)}dt.
\feq
Inspecting the coefficient of $x^ny^k/n!$ and the coefficient of $x^n/n!$ in $\witi T(x,y)$ and $\witi T(x,1)$ 
finishes the proof of Theorem~\ref{main_thm3}-(ii) and (iii). \hfill\hfill\qed
	
\subsection{Total height of weak records. Proof of Theorem~\ref{main_thm3}-(iv),(v)}
\label{wr}
Let $q_i=q^i$ for all $i\geq1$ in \eqref{Qkqs}.
Therefore,
\beq
Q_k(x;q,q^2,\ldots,q^{k-1})=x^k\frac{\prod_{j=2}^k\Big(q(1-(j-2)x)+x(q^2+\cdots+q^{j-1})\Big)}
{\prod_{j=1}^k\Big((1-x)(1-(j-1)x)-x^2(q+\cdots+q^{j-1})\Big)}.
\feq
Taking the derivative with respect to $q$ and evaluating it at $q=1,$ we obtain that
\beq
T_k(x):= \frac{\partial}{\partial q}Q_k(x;q,q^2,\ldots,q^{k-1})\,\Bigm|_{q=1}
= \frac{x^k}{\prod_{j=1}^k(1-jx)} \bigg(k-1+\binom{k}{3}+x^2\binom{k+1}{3}\bigg).
\feq
which by virtue of \eqref{SnkIdent} implies the claim in part (iv) of Theorem~\ref{main_thm3}.
\par
We turn now to the proof of part (v) of the theorem. We start with an auxiliary algebraic identity.
Recall (see, for instance, \cite{comtet} or \cite{M}) that
\beq
\frac{(e^x-1)^k}{k!}=
\sum_{n\geq k} S_{n,k}\frac{x^n}{n!}.
\feq
Using this fact we can calculate
\beq
&&\sum_{k\geq2}\sum_{n\geq k}\bigg((k-1)S_{n,k}+\binom{k}{3}S_{n-1,k}+\binom{k+1}{3}S_{n-2,k}\bigg)\frac{x^ny^k}{n!}
\\
&&
\quad
=\sum_{k\geq2}\bigg(\frac{k-1}{k!}(e^x-1)^k+\binom{k}{3}
\int_0^x\frac{(e^t-1)^k}{k!}\,dt+\binom{k+1}{3}\int_0^x\int_0^t\frac{(e^r-1)^k}{k!}\,drdt\bigg)y^k
\\
&&
\quad
=1+(y(e^x-1)-1)e^{y(e^x-1)}+\int_0^x\frac{y^3}{6}(e^t-1)^3e^{y(e^t-1)}\,dt
\\	
&&\quad \qquad +\int_0^x\int_0^t\frac{y^2}{6}(e^r-1)^2e^{y(e^r-1)}(3-y+ye^r)\,drdt.
\feq
In order to calculate the total of all heights of all weak records over $\PP_{n}$, we take advantage of the corresponding exponential generating function of $T_n(x)$ as defined by \eqref{Ae}. Summing up $\witi T_k(x)y^k$ over all $k,$ taking the second derivative of the summation with respect to $x,$ and leveraging the previous combinatorial identity, we get
\begin{align*}
\frac{\partial^2}{\partial x^2}\sum_{k\geq0}\witi T_k(x)y^k 
&=\frac{1}{6}\frac{\partial^4}{\partial x^4}e^{y(e^x-1)} +\frac{4-3y}{6}
\frac{\partial^3}{\partial x^3}e^{y(e^x-1)}+\frac{3y^2-6y-4}{6}\frac{\partial^2}{\partial x^2}e^{y(e^x-1)}
\\ 
&-\frac{y^3-3y^2+3y+1}{6}\frac{\partial}{\partial x}e^{y(e^x-1)}-\frac{y^3-3y^2}{6}e^{y(e^x-1)}.
\end{align*}
Finally, extracting the coefficient of $x^n/n!$ in $\frac{\partial^2}{\partial x^2}
\sum_{k\geq0}\witi T_k(x)$ gives Theorem~\ref{main_thm3}-(v). \hfill\hfill\qed


\begin{thebibliography}{10}	
\bibitem{Asa}
W.~Asakly, 
\emph{Sum of weighted records in set partitions}, 2019,
preprint is available at \url{https://arxiv.org/abs/1906.00680v1}.	
\filbreak

\bibitem{CMS}
N.~Cakic, T.~Mansour, and R.~Smith, 
\emph{Elements protected by records in set partitions}, 
J. Diff. Eq. Appl. \textbf{24} (2018), 1880--1893.
\filbreak 

\bibitem{Can}
E.~R.~Canfield, 
\emph{Engel's inequality for Bell numbers}, 
J. Combin. Theory. Ser. A \textbf{72} (1995), 184--187.
\filbreak 

\bibitem{hsien1}
H.-H. Chern and H.-K. Hwang,
\emph{Limit distribution of the number of consecutive records},
Random Structures Algorithms \textbf{26} (2005), 404--417.
\filbreak

\bibitem{comtet}
L.~Comtet,
\emph{Advanced Combinatorics. The Art of Finite and Infinite Expansions}, 
revised and enlarged edition, D. Reidel Publishing Co., 1974.
\filbreak

\bibitem{F}
S.~Fereti\'c, 
\emph{A perimeter enumeration of column-convex polyominoes}, 
Discrete Math. Theor. Comput. Sci. \textbf{9} (2007), 57--84.
\filbreak
	
\bibitem{flaj}
P.~Flajolet and R.~Sedgewick, 
\emph{Analytic Combinatorics}, Cambridge University Press, 2008.
\filbreak 

\bibitem{G}
A.~Geraschenko, 
\emph{An investigation of skyline polynomials}, 2016, 
preprint is available at \url{http://people.brandeis.edu/~gessel/47a/geraschenko.pdf}.
\filbreak

\bibitem{Gl}
N.~Glick, 
\emph{Breaking records and breaking boards}, 
Amer. Math. Monthly \textbf{85} (1978), 2--26.
\filbreak 
	
\bibitem{KMW}
A.~Knopfmacher, T.~Mansour, and S.~Wagner, 
\emph{Records in set partitions}, 
Electron. J. Combin. \textbf{17} (2010), \#R109.
\filbreak

\bibitem{K}
D.~E.~Knuth, 
\emph{The Art of Computer Programming. Vol. 1: Fundamental Algorithms}, 
3rd ed., Addison-Wesley, 1997.
\filbreak 

\bibitem{Ko}
I.~Kortchemski, 
\emph{Asymptotic behavior of permutation records}, 
J. Combin. Theory. Ser. A \textbf{116} (2009), 1154--1166.
\filbreak

\bibitem{Louch}
G.~Louchard,
\emph{Sum of positions of records in random permutations: asymptotic analysis}, 
Online J. Anal. Comb. \textbf{9} (2014).
\filbreak

\bibitem{M}
T.~Mansour, 
\emph{Combinatorics of Set Partitions}, 
CRC Press, 2013.
	
\bibitem{MS1}
T.~Mansour and M.~Shattuck, 
\emph{Bargraph statistics on words and set partitions}, 
J. Diff. Eq. Appl. \textbf{23} (2017), 1025--1046.
\filbreak
	
\bibitem{MW}
A.~Myers and H.~Wilf, 
\emph{Left-to-right maxima in words and multiset permutations}, 
Israel J. Math. \textbf{166} (2008), 167--183.
\filbreak 	

\bibitem{R}
A.~R\'enyi, 
\emph{Th\'eorie des \'el\'ements saillants d'une suite d'observations}, 
Ann. Fac. Sci. Univ. Clermont-Ferrand \textbf{8} (1962), 7--13.
\filbreak

\bibitem{S}
D.~Stanton and D.~White, 
\emph{Constructive Combinatorics}, Springer, 1986.
\filbreak 

\end{thebibliography}
\end{document}